\documentclass{proc-l}

\usepackage{amsmath} 
\usepackage{amsthm} 
\usepackage{amssymb} 
\usepackage{amsfonts} 
\usepackage{latexsym} 
\usepackage{graphicx} 
\usepackage{avant} 
 \usepackage[usenames,dvipsnames]{pstricks}
 \usepackage{epsfig}
 \usepackage{pst-grad} 
 \usepackage{pst-plot} 

\numberwithin{equation}{section}

\DeclareMathOperator{\vol}{vol} 
 \DeclareMathOperator{\ZZ}{\mathbb{Z}}    

 \theoremstyle{definition} 
 \theoremstyle{remark} 
 \numberwithin{equation}{section}

\title{Exact Spectral Asymptotics on the Sierpinski Gasket} 
\author{Robert S. Strichartz}
\address{Department of Mathematics, Cornell University, Ithaca, NY 14853}
\email{str@math.cornell.edu}
\thanks{Research supported in part by the National Science grant DMS-0652440}

\begin{document}

\begin{abstract}
One of the ways that analysis on fractals is more complicated than analysis on manifolds is that the asymptotic behavior of the spectral counting function $N(t)$ has a power law modulated by a nonconstant multiplicatively periodic function. Nevertheless, we show that for the Sierpinski gasket it is possible to write an exact formula, with no remainder term,  valid for almost every $t$. This is a stronger result than is valid on manifolds.

\end{abstract}

\maketitle 

\section{Introduction}

The Weyl asymptotic formula for the counting function 
\begin{equation}
	N(t)=\#\{\lambda_j\leq t\} 
\end{equation}
for the eigenvalues $\{\lambda_j\}$ of $-\Delta$ on a compact manifold $\Omega$ (with or without boundary) states that 
\begin{equation}
	N(t)\sim c_n\mbox{vol}(\Omega)t^{n/2}\;\;\;\;\;\mbox{ as }\;\;t\to\infty 
\end{equation}
where $c_n$ is a constant depending on the dimension $n$ of $\Omega$. There are more precise statements such as 
\begin{equation}
	N(t)\sim c_{n}\vol(\Omega)t^{n/2}+O(t^{(n-1)/2}) 
\end{equation}
for manifolds without boundary or with smooth boundary, and 
\begin{equation}
	N(t)=c_{n}\vol(\Omega)t^{n/2}+O(t^{d/2}) 
\end{equation}
for manifolds with non-smooth boundary of Minkowski dimension $d$ [L].

In the case of Laplacians on fractals, the situation is more complicated. For example, for the Sierpinski gasket ($SG$) with either Dirichlet or Neumann boundary conditions, the analogous statement is 
\begin{equation}
	N(t)\sim G(t)t^{\alpha}\hfill\text{ as }\hfill t\to\infty\hfill\text{ for } \alpha=\frac{\log 3}{\log 5}, 
\end{equation}
where $G$ is a multiplicatively periodic function 
\begin{equation}
	G(5t)=G(t) 
\end{equation}
that is bounded away from zero, and discontinuous (hence non-constant). This was first observed by Fukushima and Shima [FS] based on a precise description of the spectrum. It was later put into a wider context by Kigami and Lapidus [KL] based on the Reneval Theorem. Using a refinement of the Reneval Theorem, Kigami [K1] showed that the error is bounded, 
\begin{equation}
	N(t)=G(t)t^{\alpha}+O(1). 
\end{equation}
Note that this is consistent with the fact that the boundary of $SG$ consists of three points, hence has dimension zero. Similarly, the explanation for the value of $\alpha$ is that $\alpha=d/(d+1)$ where $d=\log 3/\log 5$ is the dimension of $SG$ in the resistance metric and $d+1$ is the order of the Laplacian [S2].

We already see that (1.7) is stronger than (1.3). The point of this note is that even more is true; we can have an exact formula with no remainder term at all, provided we restrict attention to almost every $t$. More precisely, let $\widetilde{SG}$ denote the double cover of $SG$ (two copies glued together at corresponding boundary points). Then $\widetilde{SG}$ is an example of a fractafold without boundary [S1]. Neumann eigenfunctions on $SG$ extended evenly to $\widetilde{SG}$ and Dirichlet eigenfunctions on $SG$ extend oddly to $\widetilde{SG}$ give a complete set of eigenfunctions on $\widetilde{SG}$. In particular 
\begin{equation}
	\widetilde{N}(t)=N_N(t)+N_D(t) 
\end{equation}
where $\widetilde{N}$, $N_{N}$ and $N_{D}$ denote the eigenvalue counting function on $\widetilde{SG}$, and $SG$ with Neumann or Dirichlet boundary conditions, respectively.

\textbf{Theorem}: \emph{There exists an open set $A\subset(0,\infty)$ whose complement has measure zero, with $5A=A$, such that for each $t\in A$ there exists $m_{0}(t)$ such that } 
\begin{equation}
	\widetilde{N}(5^{m}t)=2G(5^{m}t)(5^{m}t)^{\alpha}=2G(t)t^{\alpha}3^{m} 
\end{equation}
\emph{for all $m\geq m_{0}(t)$, and $G$ is continuous on $A$. Also}
\begin{align}
	N_{N}(5^{m}t)=G(t)t^{\alpha}3^{m}+G_{1}(t)\\
	N_{D}(5^{m}t)=G(t)t^{\alpha}3^{m}-G_{1}(t)\notag 
\end{align}
\emph{for all $m\geq m_{0}(t)$, where $G_{1}(5t)=G_{1}(t)$ and $G_{1}$ is bounded.}

Our proof is based on the ``spectral decimation'' description of the spectrum from [FS]. It seems likely that a similar result holds for many pcf fractals for which spectral decimation is valid ([Sh], [B-T]), but this would require a different proof. It is straightforward to extend the result to all fractafolds without boundary based on $SG$ [S1].

It is instructive to compare our result with the case of the second derivative on an interval. To be specific, consider the interval $[0,\pi]$ with Dirichlet boundary conditions. The eigenvalues are the perfect squares and $N(t)=[t^{1/2}]$. Then we have $N(t)=t^{1/2}$ exactly when $t$ is a perfect square, so $N(4^{m}t)=2^{m}t^{1/2}$ for all $t$ of the form $4^{-n}k^{2}$ for $m\geq n$. Thus the analog of the set $A$ is a countable dense set. This is not the first example of a result in fractal analysis that is stronger than the corresponding result in classical analysis (see [S3] for convergence of Fourier series, and also [BK], [BS], [S5]). We hope it will not be the last example.

The reader is referred to the books [K2] and [S4] for any undefined notation.

\section{Structure of the Spectrum}

We review briefly the structure of the spectrum for $\widetilde{SG}$ and $SG$ with Dirichlet or Neumann boundary conditions ([FS], [GRS], [S1], [S4]). The lowest eigenvalue is of course zero, corresponding to the constant eigenfunction, and this occurs with multiplying one for $\widetilde{SG}$ and $SG$ with Neumann boundary conditions. After that, the eigenvalues arrange themselves in cycles $C_{1}, C_{2},\ldots$. Each cycle contains ten of distinct eigenvalues. The eigenvalues and their multiplicities in a cycle are shown in Table 1. We factor out the powers of $2$ to write uniquely
\begin{equation}
k=2^{j}(2\ell-1)
\end{equation} where $j$ is a nonnegative integer and $\ell$ a positive integer. There are increasing families of primitive eigenvalues $\{\lambda^{(2)}_{n}\}$, $\{\lambda^{(3)}_{n}\}$ and $\{\lambda^{(5)}_{n}\}$ for $n=1,2,3,\ldots$. The upper index denotes the eigenvalue of the graph Laplacian of the restriction of the corresponding eigenfunctions to the nine element graph $\Gamma_{1}$ of one approximating $\widetilde{SG}$ (Figure 2.1).

\begin{figure}[htbp]
\begin{center}
\includegraphics[width=8cm]{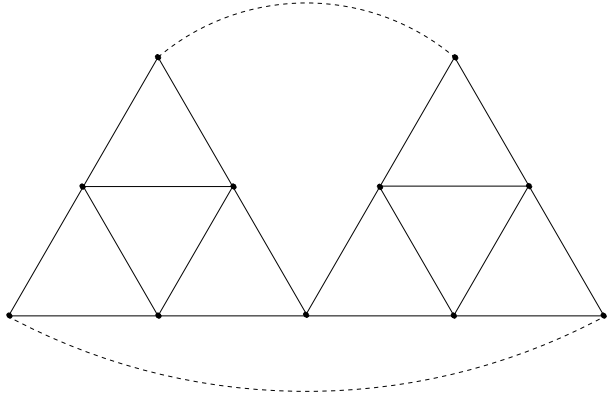}
\caption{The graph $\Gamma_{1}$ (dotted lines indicate identified points)}
\end{center}
\end{figure}

\begin{table}[htdp]
\caption{Eigenvalues and their multiplicities in the cycle $C_{k}$}
\begin{center}
\begin{tabular}{|c|c|c|c|c|}
\hline
Eigenvalues & $M_{N}$ & $M_{D}$ & $\widetilde{M}$ & $N_{N}-N_{D}$\\
\hline
$\lambda_{2k-1}^{(2)}$&0 &1 &1 &0 \\
\hline
$\lambda_{2k-1}^{(3)}$ &2&0&2&2 \\
\hline
$\lambda_{2k-1}^{(5)}$ &0&2&2&0 \\
\hline
$5\lambda_{k}^{(3)}$& 3&0&3&3\\
\hline
$\lambda_{2k}^{(5)}$&0&2&2&1\\
\hline
$\lambda_{2k}^{(3)}$&2&0&2&3\\
\hline
$\lambda_{2k}^{(2)}$&0&1&1&2\\
\hline
$5^{j+1}\lambda_{2\ell-1}^{(5)}$& $\frac{3^{j+1}-1}{2}$ & $\frac{3^{j+1}+3}{2}$&$3^{j+1}+1$ & 0\\
\hline
$5^{j+2}\lambda_{\ell}^{(3)}$ & $\frac{3^{j+2}+3}{2}$ & $\frac{3^{j+2}-3}{2}$&$3^{j+2}$&3 \\
\hline
$5^{j+1}\lambda_{2\ell-1}^{(5)}$& $\frac{3^{j+1}-1}{2}$ & $\frac{3^{j+1}+3}{2}$&$3^{j+1}+1$ & 1\\
\hline
\end{tabular}
\end{center}
\label{default}
\end{table}

The table lists the eigenvalues in increasing order (top to bottom), followed by their multiplicities in each of the three cases. Note that $\widetilde{M}=M_{N}+M_{D}$. The last column gives the difference of the eigenvalue counting functions $N_{N}(t)-N_{D}(t)$ for $t$ equal to the eigenvalue. Since this is a cumulative total, it needs to be justified by an induction argument. Note that $N_{N}(0)-N_{D}(0)=1$, so the values are correct for the first cycle. Since the value is 1 at the end of the first cycle, the pattern repeats at each subsequent cycle. Note that we know a priori that $0\leq N_{N}(t)-N_{D}(t)\leq 3$ because the boundary of $SG$ consists of three points,, and this is consistent with the values in the Table.

Note that the multiplicities are the same for the first seven eigenvalues of each cycle, but the last three entries in the Table have multiplicities that depend on the value of $j$ in (2.1).

We can also say explicitly what the values of the values of the primitive eigenvalues are. Define

\begin{align}
\varphi_{\pm}(t)&=\frac{5\pm\sqrt{25-4t}}{2}\;\;\;\;\;\;\;\text{ and }\\
\psi(t)&=\frac{3}{2}\lim_{k\to\infty}5^{k}\varphi_{-}^{(k)}(t)
\end{align} where $\varphi_{-}^{(k)}$ denotes $\varphi_{-}$ composed with itself $k$ times (the convergence of the limit in (2.3) follows from the Taylor expansion $\varphi_{-}(t)=\frac{1}{5}t+O(t^{2})$ near $t=0$). More generally, if $\delta=(\delta_{1},\delta_{2},\ldots,\delta_{n})$ is a finite sequence of $\pm$ signs, we write $\varphi_{\delta}=\varphi_{\delta_{1}}\circ\varphi_{\delta_{2}}\circ\cdots\circ \varphi_{\delta_{m}}$ and $|\delta|=m$. Then
\begin{equation}
\{\lambda_{n}^{(p)}\}=\{5\psi(p),5^{2}\psi(\varphi_{+}(p)),5^{3}\psi(\varphi_{++}(p)), 5^{3}\psi(\varphi_{+-}(p))),\ldots,5^{|\delta|+1}\psi(\varphi_{\delta}(p)),\ldots\}
\end{equation} for $p=2,3$ or 5, and $\delta_{1}=+$ but otherwise $\delta$ is unrestricted. Therefore there are $2^{m-1}$ choices of $\delta$ with $|\delta|=m$, and since $\delta_{-}$ is increasing and $\delta_{+}$ is decreasing it is necessary to use an unconventional ordering in order for the list (2.4) to be in increasing order. However we always have $5^{|\delta|+1}\psi(\varphi_{\delta}(p))<5^{|\delta'|+1}\psi(\varphi_{\delta'}(p'))$ if $|\delta|<|\delta'|$. This follows from the estimate 
\begin{equation}
\psi(3)5^{|\delta|+1}\leq 5^{|\delta|+1}\psi(\varphi_{\delta}(p))< \psi(5)5^{|\delta|+1}
\end{equation} which follows from the estimate $3\leq\varphi_{+}(t)\leq 5$ on $0\leq t\leq 6$ and the fact that $\psi$ is increasing. 

The following is the key Lemma in the proof of the Theorem. Note that it is really (1.9) for the special values $t=\lambda^{(3)}_{n}$.

\textbf{Lemma:} For any $n$ and $j$,
\begin{equation}
\widetilde{N}(5^{j+2}\lambda_{n}^{(3)})=3^{j}\widetilde{N}(5^{2}\lambda^{(3)}_{n}).
\end{equation}
\begin{proof} For $n=1$ it is known that $\widetilde{N}(5^{j+2}\lambda_{1}^{(3)})=3^{j+3}$, the total number of vertices in the level $j+2$ graph approximation $\Gamma_{j+2}$ of $\widetilde{SG}$, It is easy to see in any case from the Table that $\widetilde{N}(5^{2}\lambda^{(3)}_{1})=27$ and $\widetilde{N}(5^{3}\lambda_{1}^{(3)})=81.$ We will show by induction on $k$ that
\begin{equation}
\widetilde{N}(5^{j+3}\lambda_{\ell}^{(3)})=3\widetilde{N}(5^{j+2}\lambda_{\ell}^{(3)})
\end{equation} and this verifies the $k=1$ case. So we assume that the result is true for $k-1$. The argument is slightly different depending on the parity of $k$, so assume first that $k=2\ell-1$ is odd. Then $k-1=2^{j'}(2\ell'-1)$ is even, and the induction hypothesis is
\begin{equation}
\widetilde{N}(5^{j'+3}\lambda_{\ell'}^{(3)})=3\widetilde{N}(5^{j'+2}\lambda_{\ell'}^{(3)}).
\end{equation} What we need to show is 
\begin{equation}
\widetilde{N}(5^{3}\lambda_{\ell}^{(3)})-\widetilde{N}(5^{j'+3}\lambda_{\ell'}^{(3)})=3\left(\widetilde{N}(5^{j'}\lambda_{\ell}^{(3)})-\widetilde{N}(5^{j'+2}\lambda_{\ell'}^{(3)})\right).
\end{equation} Note that on the right side of (2.9) we are counting the last eigenvalue in cycle $C_{k-1}$ (multiplicity $3^{j'+1}+1$) and the first nine eigenvalues of cycle $C_{k}$ (multiplicities $1+2+2+3+2+2+1+4+9=26$), so the right side of (2.9) is $3^{j'+2}+3^{4}$. On the left side of (2.9) we are counting the last eigenvalue in cycle $C_{2k-2}$ (multiplicity $3^{j'+2}+1$), all the eigenvalues in cycle $C_{2k-1}$ (multiplicities $1+2+2+3+2+2+1+4+9+4=30$) and the first nine eigenvalues in cycle $C_{2k}$ (multiplicities $1+2+2+3+2+2+1+10+27=50$) for a total of $81+3^{j'+2}$. This verifies (2.9).

Finally, assume $k=2^{j}(2\ell-1)$ is even, and $k-1=2\ell'-1$ is odd. The induction hypothesis is 
\begin{equation}
\widetilde{N}(5^{3}\lambda_{\ell'}^{(3)})=3\widetilde{N}(5^{2}\lambda_{\ell'}^{(3)}), \text { and we need to show }
\end{equation}
\begin{equation}
\widetilde{N}(5^{j+3}\lambda_{\ell}^{(3)})-\widetilde{N}(5^{3}\lambda_{\ell'}^{(3)})=3\left(\widetilde{N}(5^{j+2}\lambda_{\ell}^{(3)}-\widetilde{N}(5^{2}\lambda_{\ell'}^{(3)})\right).
\end{equation} We are counting the same eigenvalues as before but the multiplicities are different. On the right side the multiplicities are 4 and $1+2+2+3+2+2+1+(3^{j+1}+1)+3^{j+2}=14+3^{j+1}+3^{j+2}$. On the left side the multiplicities are 10 from cycle $C_{2k-2}$, $1+2+2+3+2+2+1+4+9+4=30$ from cycle $C_{2k-1}$, and $1+2+2+3+2+2+1+(3^{j+2}+1)+3^{j+3}=14+3^{j+2}+3^{j+3}$ from cycle $C_{2k}$, for a grand total of $54+3^{j+2}+3^{j+3}$.
\end{proof}

There appears to be a lot of arithmetic in this proof with no clear explanation. Certainly it would be desirable to have a more conceptual proof. 

\section{Proof of Theorem}

We consider two intervals between $5^{2}\lambda_{\ell}^{(3)}$ and its neighboring eigenvalues, namely
\begin{equation}
\begin{cases}
A_{\ell}= \left( 5\lambda_{2\ell-1}^{(5)}\;,\;5^{2}\lambda_{\ell}^{(3)} \right)  \\
A'_{\ell}=\left( 5^{2}\lambda_{\ell}^{(3)}\;,\;5\lambda_{2\ell}^{(5)} \right)
\end{cases}
\end{equation} The function $\widetilde{N}(t)$ is constant on these intervals, 
\begin{equation}
\widetilde{N}(t)=
\begin{cases}
 \widetilde{N}(5^{2}\lambda_{\ell}^{(3)})-9,  & t\in A_{\ell}\\
 \widetilde{N}(5^{2}\lambda_{\ell}^{(3)}),  & t\in A'_{\ell}.
\end{cases}
\end{equation} More generally, $\widetilde{N}(t)$ is constant on the intervals $5^{j}A_{\ell}$ and $5^{j}A'_{\ell}$ for $j\geq0$, namely 
\begin{equation}
\widetilde{N}(5^{j}t)=
\begin{cases}
 3^{j}\widetilde{N}(5^{2}\lambda_{\ell}^{(3)})-3^{j+2},  & t\in A_{\ell}\\
3^{j} \widetilde{N}(5^{2}\lambda_{\ell}^{(3)}),  & t\in A'_{\ell}.
\end{cases}
\end{equation}This uses the Lemma and the fact that the multiplicity of $5^{j+2}\lambda_{\ell}^{(3)}$ is $3^{j+2}$. This gives us (1.9) with $m_{0}(t)=0$ on $A_{\ell}\cup A'_{\ell}$ with 
\begin{equation}
G(t)=
\begin{cases}
 \frac{1}{2}\left(\widetilde{N}(5^{2}\lambda_{\ell}^{(3)})-9\right)t^{-\alpha},  & t\in A_{\ell}\\
\frac{1}{2}\widetilde{N}(5^{2}\lambda_{\ell}^{(3)})t^{-\alpha},  & t\in A'_{\ell}.
\end{cases}
\end{equation}For any $n\in\ZZ$, we similarly obtain (1.9) with $m_{0}(t)=\max\{-n,0\}$ for $t\in 5^{n}(A_{\ell}\cup A_{\ell'})$ with $G$ extended by multiplicative periodicity. We also observe from the Table that (1.10) holds on the same intervals with
\begin{equation}
G_{1}(t)=
\begin{cases}
0,  & t\in A_{\ell}\\
\frac{3}{2}, & t\in A'_{\ell}.
\end{cases}
\end{equation}

Now let $A=\bigcup_{n=-\infty}^{\infty}\bigcup_{\ell=1}^{\infty}5^{n}\left(A_{\ell}\cup A'_{\ell}\right)$. To complete the proof we have to show that the complement of $A$ has measure zero. We can restrict attention to the interval $[0,5]$ because $5A=A$. Write $A''_{\ell}=A_{\ell}\cup\{5^{2}\lambda^{(3)}_{\ell}\}\cup A'_{\ell}=\left(5^{(5)}\lambda_{2\ell-1}^{(5)}\;,\;5\lambda_{2\ell}^{(5)} \right)$ and $A''=\bigcup_{n=-\infty}^{\infty}\bigcup_{\ell=1}^{\infty}5^{n}A''_{\ell}$. The complement of $A''$ differs from the complement of $A$ by a countable set of points, and it simplifies matters to deal with $A''$. Now each of the intervals $A_{\ell}''$ has the form 
\[(5^{m}\psi(\varphi_{\delta}(5)),5^{m}\psi(\varphi_{\delta'}(5))\]
where $|\delta|=|\delta'|=m-1$ and we no longer require that $\delta_{1}$ and $\delta'_{1}$ be $+$ (this uses the fact that $\psi(\varphi_{-}(x))=5\psi(\varphi(x)$). Then $5^{-m}A''_{\ell}=(\psi(\varphi_{\delta}(5),\psi(\varphi_{\delta'}(5))\subseteq[0,5]$. The composition with $\psi$ is just a differentiable distortion factor, so we drop it and define $B_{\ell}=(\varphi_{\delta}(5),\varphi_{\delta'}(5))$ and $B=\bigcup_{\ell=1}^{\infty}B_{\ell}$. Our problem is to show that the complement of $B$ in $[0,5]$ has measure zero.

Call this complement $J$. We claim that $J$ is just the Cantor set defined by the mappings $\varphi_{+}$ and $\varphi_{-}$ on $[0,5]$, or equivalently the Julia set for the mapping $z\to z(5-z)$ since $\varphi_{\pm}$ are the the inverses of this mapping. It is known that this Julia set has measure zero, but in any case this is easily seen from estimates on the derivatives of $\varphi_{\pm}$.

Note that $B_{1}=\left(\frac{5-\sqrt{5}}{2},\frac{5+\sqrt{5}}{2}\right)$ so the complement of $B_{1}$ is $\varphi_{-}([0,5])\cup\varphi_{+}([0,5])$. Then $B_{2}=\varphi_{-}(B_{1})$ and $B_{3}=\varphi_{+}(B_{1})$ so the complement of $B_{1}\cup B_{2} \cup B_{3}$ is $\varphi_{-  -}([0,5])\cup\varphi_{-+}([0,5])\cup\varphi_{++}([0,5])\cup\varphi_{+-}([0,5])$. Iterating this reasoning, the complement of $\bigcup_{\ell=1}^{2^{m}-1}B_{\ell}$ is $\bigcup_{|\delta|=m}\varphi_{\delta}([0,5])$. This verifies the Cantor set description of $J$.\hfill$\square$

\end{document}